\documentstyle[11pt]{article}
\newtheorem {th}{Theorem}[section]
\newtheorem {lem}[th]{Lemma}
\newtheorem {pr}[th]{Proposition}

\def\Cox{\hfill \Box}

\def\sign{\mbox{sgn}}
\def\grad{\bigtriangledown}
\def\deq{\, {\stackrel {def} {=}}}

\def\sf{\sigma\mbox{-field}}
\def\ee{\epsilon}
\def\dd{\delta}

\def\cL{{\cal L}}
\def\LL{{\cal L}}
\def\S{{\cal S}}
\def\C{{\cal C}}
\def\B{{\cal B}}
\def\E{{\bf{E}}}
\def\P{{\bf{P}}}
\def\Q{{\bf{Q}}}

\def\R{{\bf{R}}}
\def\Z{{\bf{Z}}}
\def\F{{\cal{F}}}

\def\|{\, | \, }
\def\one{{\bf 1}}

\begin{document}

\begin{center}
{\Large \bf The probability that Brownian motion almost contains a line}
\end{center}

\begin{flushright}
Robin Pemantle\footnote{Research supported in part by a grant from the
Alfred P.\ Sloan Foundation, by a Presidential Faculty Fellowship and by
NSF grant \# DMS9300191} \\
Department of Mathematics \\
University of Wisconsin-Madison \\
pemantle@math.wisc.edu
\end{flushright}
\vspace{.2in}

{\bf ABSTRACT:}

\noindent{Let} $R_\ee$ be the Wiener sausage of radius $\ee$ 
about a planar Brownian motion run to time 1.  Let $P_\ee$ be
the probability that $R_\ee$ contains the line segment $[0,1]$ 
on the $x$-axis.
Upper and lower estimates are given for $P_\ee$.  The upper 
estimate implies that the range of a planar Brownian motion 
contains no line segment.

\section{Introduction}

Let $\{ B_t : t \geq 0 \} = \{ (X_t , Y_t) \}$ be a $2$-dimensional 
Brownian motion, started from the origin unless otherwise stated. 
Let $R = B[0,1]$ be the range of this trajectory up to time 1.  
Brownian motion is one of the most fundamental stochastic processes, 
applicable to a wide variety of physical models, and consequently its 
properties have been widely studied.  Some people prefer to
study simple random walk, but questions in the discrete case are
nearly always convertible to the continuous case: rescale the
simple random walk $\{ X_n : n \in \Z \}$ in space and time to
give a trajectory $\{ N^{-1/2} X_n : n \in N^{-1} \Z \}$ and 
replace the question of whether this trajectory covers or avoids
a lattice set $A \subseteq (N^{-1} \Z)^2$ by the question of
whether the Wiener sausage $\{ B_t + x : 0 \leq t \leq 1 , 
|x| \leq N^{-1/2} \}$ covers or avoids the set $A + [-N^{-1/2} ,
N^{-1/2}]^2$.  

Focusing on the random
set $R$ is a way of concentrating on the geometric properties
of Brownian motion and ignoring those having to do with time.
The sorts of things we know about $R$ are as follows.  Its
Hausdorff dimension is 2, and we know the exact gauge function
for which the Hausdorff measure is almost surely positive and 
finite \cite{Ra}.  Studies from 1940 to the present of intersections 
of independent Brownian paths (or equivalently multiple points
of a single path)
show that two paths are always disjoint in dimensions 4 and higher,
two paths intersect in three dimensions but three paths do not, 
while in 2 dimensions any finite number of independent copies of $R$ 
have a common intersection with positive probability~\cite{Lev,DEK,LeG}.

Questions about which sets $R$ intersects are reasonably 
well understood via potential theory.  For example, Kakutani
showed in 1944 that in any dimension
$$\P [R \cap S \neq \emptyset] > 0 \Leftrightarrow cap_G (S) > 0 ,$$
where $cap_G$ is the capacity with respect to the Green's kernel
\cite{Kak}.  There are versions of this result for intersections
with several independent copies of $R$ \cite{FS,Sa} and quantitative
estimates exist (up to computable constant factors) for  
$\P [R \cap S \neq \emptyset]$ \cite {BPP}.  Not so well 
understood are the dual questions to these, namely covering
probabilities.  The most interesting questions seem to arise in
dimension 2.  In the discrete setting, 
one may ask for the radius $L$ of the largest disk about the origin
covered by a simple random walk up to time $N$.  The answer is
roughly $\exp (c \sqrt{\log N})$, in the sense that 
$\P [\ln L \geq t \sqrt{\log N}]$ is bounded away from 0 and 1
for fixed $t$ as $N \rightarrow \infty$; see~\cite{La,Re,Anon}.
Another type of covering
question is whether $R$ contains any set from a certain class.
For example, it is known that $R$ contains a self-avoiding
path of Hausdorff dimension 2 \cite{Bu3} and a self-avoiding
path of dimension less than $3/2$ \cite{BL}, but it is
not known whether the infimum of dimensions of self-avoiding paths
contained in $R$ is 1.  A modest start on an answer would be
to show something that seems obvious, namely that $R$ almost 
surely contains no line segment.  

Despite the seeming self-evidence of this assertion, I know of 
no better proof than via the following Wiener sausage estimates.
Let 
$$R_\ee = \{ x \in \R^2 : |y - x| \leq \ee \mbox{ for some } 
   y \in R \}$$
be the Wiener sausage of radius $\ee$.  Let 
$$\P_\ee = \P [R_\ee \supseteq [0,1] \times \{ 0 \} ] $$ 
be the probability that the Wiener sausage contains the unit 
interval on the $x$-axis.  An easy upper bound for $\P_\ee$
is given by $\P [ (1/2 , 0) \in R_\ee]$ which is order
$1 / |\log \ee|$.  Omer Adelman (personal communication)
has obtained an upper bound that is $o(\ee^2)$, though this
was never written down; see Remark~2 following Theorem~\ref{th 2}.
An obvious lower bound of order $\exp (-c / \ee)$ is gotten
from large deviation theory by forcing the Brownian motion 
to hit each of $\ee^{-1}$ small disks in order.  The main 
result of this paper is to improve the upper and lower 
bounds on $\P_\ee$ so as to be substantially closer to each other.
\begin{th} \label{th 1}
There exist positive constants $c_1, c_2, c_3$ and $c_4$ such
that for all sufficiently small $\ee$ the following four 
inequalities hold:
\begin{equation} \label{eq upper}
\P_\ee \leq c_1 \exp ( - |\log \ee|^2 / c_2 \log^2 |\log \ee|)
\end{equation}
and
\begin{equation} \label{eq lower}
\exp ( - c_4 |\log \ee|^4 ) \leq \P_\ee ;
\end{equation}
if $\Xi$ is the Lebesgue measure of the 
intersection of $R_\ee$ with any fixed line segment of length at 
most 1 contained in $[0,1]^2$, then for any $\theta \in (0,1)$, 
\begin{equation} \label{eq upmeasure}
   \P [\Xi \geq \theta] \leq \exp (- |\log \ee|^2 \theta^2 / c_3 
   \log^2 |\log \ee|)
\end{equation}
and
\begin{equation} \label{eq lomeasure}
\exp (- c_4 |\log {1 - \theta \over 3}|^2 |\log \ee |^2) \leq
   \P [\Xi \geq \theta] .
\end{equation}
\end{th}
This implies
\begin{th} \label{th 2}
A 2-dimensional Brownian motion, run for infinite time,
almost surely contains no line segment.  In fact it intersects 
no line in a set of positive Lebesgue measure.  
\end{th}

{\em Remarks:}

\noindent{1. } A discrete version of these results holds,
namely that simple random walk run for time $N^2$ covers
the set $\{ (1 , 0) , \ldots , (N , 0) \}$ with probability
between $\exp (- C_4 (\log N)^4)$ and $C_1 \exp(
- \log^2 N / C_2 log^2 (\log N))$; similarly the probability
of covering at least $\theta N$ of these $N$ points is between
$$\exp \left ( - C_4 \left | \log {1 - \theta \over 3} \right |^2 
   \log^2 N \right )$$
and
$$\exp (- \log^2 N \theta^2 / C_3 \log^2 \log N) \, .$$
The proof is entirely analogous (though the details are a little 
trickier) and is available from the author.

\noindent{2. } The proof of Theorem~\ref{th 2} requires only
the upper bound $\P_\ee = o (\ee^2)$, which is weaker
than what is provided by Theorem~\ref{th 1}.  Thus Adelman's
unpublished bound, if correct, would also imply Theorem~\ref{th 2}.
However, since both Adelman's proof and the proof of 
Theorem~\ref{th 1} are non-elementary (Adelman discretizes 
as in Remark~1 and then sums over all possible orders in
which to visit the $N$ points), the question of Y.\ Peres which
motivated the present article is still open: can you find an 
elementary argument to show $\P_\ee = o(\ee^2)$?  

This section concludes with a discussion of why 
Theorem~\ref{th 2} follows
from the upper estimate on $\P_\ee$.  The next section begins
the proof of the upper bounds~(\ref{eq upper}) 
and~(\ref{eq upmeasure}), first replacing time 1 by a hitting time
and then proving some technical but routine lemmas on 
the quadratic variation of a mixture of stopped martingales.
Section~3 finishes the proof of the upper bounds and Section~4
proves the lower bounds~(\ref{eq lower}) and~(\ref{eq lomeasure}).

\noindent{\sc Proof of Theorem~2 from Theorem~1}:  It suffices
by countable additivity to show that $R \cap [0,1]^2$ 
almost surely intersects no line in positive Lebesgue measure. 
Let $m(\cdot)$ denote 1-dimensional Lebesgue measure; then
it suffices to see for each $\theta > 0$ that $\P [H_\theta] = 0$, 
where $H_\theta$ is the event that for some line $l$,
$m(R \cap [0,1]^2 \cap l) \geq \theta$.  From now on, fix $\theta > 0$.

Let $\cL_\ee$ be the set of line segments whose endpoints are 
on the boundary of $[0,1]^2$ and have coordinates that are 
integer multiples of $\ee$.  The set $\cL_\ee'$ of halves of line segments
in $\cL_\ee$ has twice the cardinality and contains segments
of maximum length $\sqrt{2} / 2$.  For any line $l$ intersecting
the interior of $[0,1]^2$ and any $\ee$ there is a line segment 
$l_\ee \in \cL_\ee$ whose endpoints are within $\ee$ of the 
respective two points where $l$ intersects $\partial [0,1]^2$.  
If $m (R \cap [0,1]^2 \cap l) \geq \theta$ then $m(R_\ee \cap l_\ee)$ 
will be at least $\theta - 2 \ee$.  Thus on the event $H_\theta$,
every $\ee > 0$ satisfies 
$$m (R_\ee \cap l) \geq \theta - 2 \ee \mbox{ for some } l \in \cL_\ee ,$$
and hence 
$$m (R_\ee \cap l) \geq \theta / 2 - \ee \mbox{ for some } l \in \cL_\ee' .$$
By~(\ref{eq upmeasure}), when $\ee < \theta / 3$, the probability of 
this event is at most 
$$2|\cL_\ee | \exp (- |\log \ee|^2 (\theta / 6)^2 / c_3
   \log^2 |\log \ee|) .$$
Since $|\cL _\ee|$ is of order $\ee^{-2}$, this goes to zero
as $\ee \rightarrow 0$ which implies $\P [H_\theta] = 0$,
proving the theorem.   $\Cox$

\section{Change of stopping time and martingale lemmas}

The proof of the upper bound~(\ref{eq upmeasure}), which implies
the other upper bound~(\ref{eq upper}), comes in three pieces.  
It is convenient to stop the process
at the first time the $y$-coordinate of the Brownian motion 
leaves the interval $(-1,1)$.  The first step is to justify
such a reduction.  A reduction is made at the same time to
the case where $l$ is the line segment $[0,1] \times \{ 0 \}$.
The proof then proceeds by representing the
measure $Z := m(R_\ee \cap l)$ as the final term in a martingale
$\{ Z_t : 0 \leq t \leq 1 \}$, where $Z_t$ is simply 
$\E (Z \| \F_t)$ and $\F_t$ is the natural filtration
of the Brownian motion.  The second step is to identify the
quadratic variation of this martingale.  This involves 
a few technical but essentially routine lemmas.  The final
step, carried out in Section~3, is to obtain upper bounds for 
the quadratic variation which lead directly to tail bounds 
on $Z - \E Z$.  

Let $l_0$ be the line segment $[0,1] \times \{ 0 \}$; here is
why it suffices to prove~(\ref{eq upmeasure}) for the line $l_0$.
Let $l$ be any line of length at most 1 contained in $[0,1]^2$.
Let $\sigma$ be the first time $l$ is hit.  Let $l_1$ and $l_2$
be the two closed line segments whose union is $l$ and
whose intersection is $B_\sigma$.  If the measure of 
$R_\ee \cap l$ is at least $\theta$ then the measure of
$R_\ee \cap l_j$ is at least $\theta / 2$ for some $j$, and
the probability of this is at most twice the probability
of $R_\ee \cap l_0$ having measure at least $\theta / 2$.
Thus if $c_3$ works in~(\ref{eq upmeasure}) for $l_0$,
then $5 c_3$ works for arbitrary $l$.  

The more serious reduction is to define a stopping time
$$\tau = \inf \{ t : |Y_t| \geq 1 \}$$
and to replace the event $\{ \Xi \geq \theta \}$ with the
event $\{ \Xi \geq \theta \} \cap \{ \tau \geq 1 \}$.  The
value of this replacement will be seen in the next section,
where some estimates depend on stopping at a hitting time.
To carry out this reduction, let $G_1$ be the event
$\{ \tau \geq 1 \}$, let $K = \P [G_1]^{-1}$, and let
$G_2^{(\theta)}$ be the event $\{ \Xi \geq \theta \}$.  
The requisite lemma is:
\begin{lem} \label{lem strip}
$$\P [G_2^{(\theta)}] \leq K \P [G_1 \cap G_2^{(\theta)}] .$$
\end{lem}

\noindent{\sc Proof:}  I must show that
\begin{equation} \label{eq 3}
\P [G_2^{(\theta)} \| G_1] \geq \P [G_2^{(\theta)}] .
\end{equation}
The law of Brownian motion up to time 1, conditioned on the event
$G_1$ is a time inhomogenous Markov process (an $h$-process) and
is representable as a Brownian motion with drift.  By translational
symmetry in the $x$-direction, the drift, $(h_1 (x , y , t) , 
h_2 (x , y , t))$ must be toward the $x$-axis, i.e.,
$h_1 = 0$ and $y h_2 \leq 0$.  Given an unconditioned
Brownian motion, $\{ (X_t , Y_t) \}$ with $Y_0 \in [0,1]$,
we may in fact construct this $h$-process $\{ W_t = (X_t , V_t) \}$ 
by defining $V_0 = Y_0$ and define $\{ V_t \}$ by
$$dV_t = \sign(V_t) dY_t + h(Y_t , t) dt .$$
Then $(X,V)$ is distributed as a Brownian motion conditioned 
to stay inside the strip until time 1 and is coupled so that
$|V_t| \leq |Y_t|$
for all $t \leq \tau \wedge 1$.  The coupling is constructed
such that for all $t \leq 1$ and $\ee > 0$,
$$\{ x : |B_s - (x,0)| < \ee \mbox{ for some } s \leq t \} \subseteq
  \{ x : |W_s - (x,0)| < \ee \mbox{ for some } s \leq t \} .$$
Equation~(\ref{eq 3}) follows directly.   $\Cox$

Let, with the natural identification of $l_0$ and $[0,1]$, 
$$A_t = B[0 , t \wedge \tau ]_\ee \cap l_0 = \{ x \in [0,1] : 
   |B_s - (x,0)| \leq \ee \mbox{ for some } s \leq t \wedge \tau \} .$$
For achieving stochastic upper bounds on $m(A_1)$, one may take
advantage of the trivial inequality $m(A_1) \leq m(A_\infty)$.  Thus
let $M_t = \E (m(A_\infty) \| \F_t)$.  On the
event $G_1 \cap G_2^{(\theta)}$, clearly $M_1 \geq \theta$.  By the
previous lemma, the inequality~(\ref{eq upmeasure}) will
follow from the inequality
\begin{equation} \label{eq to show}
\P [M_1 \geq \theta ] \leq \exp ( -|\log \ee|^2 \theta^2
   / c_3' \log^2 |\log \ee|) 
\end{equation}
with $c_3 > c_3'$ and $\ee$ small enough to make up for the
factor of $K$.  

Now comes the second step, namely identifying the quadratic 
variation of the martingale $M$.  This step is just the
rigorizing (Lemma~\ref{lem integrate}) of the following intuitive 
description of $M$.
\begin{quotation}
If $I(x,t,\ee)$ is the indicator function of the event that the 
Brownian motion, killed at time $\tau$, has visited the $\ee$-ball
around $x$ by time $t$, then by the Markov property, 
if $\tau$ has not yet been reached,
\begin{eqnarray*}
M_t & = & m (A_t) + \int_{{[0,1]} \setminus A_t} \P_{B_t} 
   [x \in A_\infty] \, dx \\[1ex]
& = & m(A_t) + \int_{[0,1]} [1 - I (x, t ,\ee)] 
   \E_{B_t} I (x , \infty , \ee) \, dx . 
\end{eqnarray*}
Thus the quadratic variation $d <M>_t$ is given by 
$$ dt \left | \grad \int_{x \in [0,1] \setminus A_t} g_\ee 
   (B_t , x) dx \right |^2 ,$$
where 
$$g_\ee (y,x) = \E_y [I(x, \infty ,\ee)] $$
and $\grad$ denotes the gradient with respect to the $y$.
\end{quotation}
The formal statement of this is: 

\begin{lem} \label{lem integrate}
Let $D$ be the closure of a connected open set $D^o \subset \R^d$
and let $(\S, \B, \mu)$ be a probability space.  Suppose that $\partial D$
is nonpolar for Brownian motion started inside $D$ and that $f : D 
\times \S \rightarrow [0,1]$ has the following properties:
\begin{quote}
$(i)$ $f(x , \alpha)$ is jointly measurable, continuous in $x$ for
each $\alpha$, and $f(\cdot , \alpha)$ is identically zero on
$\partial D$;

$(ii)$ There is a constant $C$ such that $f( \cdot , \alpha)$ is
Lipschitz on $D$ with constant $C$ for all $\alpha$;

$(iii)$ For each $\alpha$, let $D_\alpha = \{ x \in D : f (x,\alpha)
\in (0,1) \}$; then for each $\alpha$ $f( \cdot , \alpha)$ is 
$C^2$ on $D_\alpha$ and its Laplacian vanishes there.  
\end{quote}
Let $\{ B_t \}$ be a Brownian motion started from some $B_0 \in D^o$.
For each $\alpha$, let $\tau_\alpha = \inf \{ t : f(B_t , \alpha)
\in \{ 0 , 1 \} \}$ be the first time $B_t$ exits $D_\alpha$ 
(in particular, $\tau_\alpha$ is at most $\tau_{\partial D}$, 
the exit time of $D^o$).  Then
$$M_t := \int f (B_{t \wedge \tau_\alpha} , \alpha) \, d\mu (\alpha)$$
is a continuous $\{ \F_t \}$-martingale with quadratic variation
$$<M>_t = \int_0^t \left | \int_\S \one_{s < \tau_\alpha} 
   \grad f (B_s , \alpha) \, d\mu (\alpha) \right |^2 \, ds .$$
\end{lem}

The proof of the lemma requires two elementary 
propositions, whose proofs are routine and relegated to the appendix.

\begin{pr} \label{pr COV}
With the notation and hypotheses of Lemma~\ref{lem integrate}, 
for each $\alpha \in \S$
let $M^{(\alpha)}_t = f (B_{t \wedge \tau_\alpha} , \alpha)$.  
Then $M^{(\alpha)}$ is a continuous $\{ \F_t \}$-martingale 
for each $\alpha$ and the brackets satisfy
$$[M^{(\alpha)} , M^{(\beta)}]_t = \int_0^t \one_{s < \tau_\alpha}
   \one_{s < \tau_\beta} \grad f (B_s , \alpha)
   \cdot \grad f (B_s , \beta)  \, ds . $$
$\Cox$
\end{pr}

\begin{pr} \label{pr together}
Let $\{ M^{(\alpha)}_t \}$ be continuous $\{ \F_t \}$-martingales as 
$\alpha$ ranges over the non-atomic measure space $(\S , \B , \mu)$, 
jointly measurable in $t$ and $\alpha$ and taking values in $[0,1]$.
Suppose that the infinitesimal covariances are all bounded,
that is, for some $C'$, all $\alpha, \beta \in \S$ and all $t > s$, 
$$\E \left [ ( M^{(\alpha)}_t - M^{(\alpha)}_s) 
   ( M^{(\beta)}_t - M^{(\beta)}_s) \| \F_s \right ] \leq C' (t-s) . $$
Then $M := \int M^{(\alpha)}
\, d\mu (\alpha)$ is a continuous martingale with
$$<M>_t = \int_0^t \left [ {\int \int}_{\S^2} d[M^{(\alpha)} , M^{(\beta)} 
  ]_s \, d\mu(\alpha) \, d\mu (\beta) \right ] .$$
$\Cox$
\end{pr}

\noindent{\sc Proof of Lemma}~\ref{lem integrate}:  Applying
Propositions~\ref{pr COV} and~\ref{pr together} to the collection
$\{ M^{(\alpha)}_t \}$ shows that $M$ is a continuous martingale
with covariance
\begin{equation} \label{eq 1}
<M>_t = \int_0^t \left [ {\int \int}_{\S^2} \one_{s < \tau_\alpha}
   \one_{s < \tau_\beta} \grad f (B_s , \alpha) \cdot \grad f 
   (B_s , \beta) \, d\mu (\alpha) \, d\mu (\beta) \right ]  \, ds . 
\end{equation}
If $s < \tau_\alpha , \tau_\beta$ then $B_s \in D_\alpha \cap D_\beta$,
so $\grad f (B_s , \alpha)$ and $\grad f (B_s , \beta)$ exist and
have modulus at most $C$ by the Lipschitz condition.  Thus
$${\int \int}_{\S^2} \one_{s < \tau_\alpha} \one_{s < \tau_\beta}
   |\grad f (B_s , \alpha)| |\grad f (B_s , \beta)| 
   \, d\mu (\alpha) \, d \mu (\beta) < \infty $$
which implies absolute convergence of the inner integral 
in~(\ref{eq 1}) and thus its factorization as
$$ {\int \int}_{\S^2} \one_{s < \tau_\alpha}
   \one_{s < \tau_\beta} \grad f (B_s , \alpha) \cdot \grad f 
   (B_s , \beta) \, d\mu (\alpha) \, d\mu (\beta)  =
   \left | \int_\S \one_{s < \tau_\alpha} \grad f (B_s , \alpha)
   \, d\mu (\alpha) \right |^2 ,$$
which proves the lemma.   $\Cox$

\section{Proof of upper bounds}

Let $D \subseteq \R^2$ be the strip $\{ (x_1 , x_2) : |x_2| 
\leq 1 \}$.  Let $\{ B_t \}$ be a 2-dimensional Brownian
motion, with $B_t = (X_t , Y_t)$ as before, and let
$\tau$ be the hitting time of $\partial D$.  Recall the notation
$$A_t = \{ x \in l_0 : |B_s - (x,0)| \leq \ee \mbox{ for some }
   s \leq t \wedge \tau \} .$$
Fix $\ee > 0$ and for $\alpha \in \R$ define
$$f (x , \alpha) = \P_x [ (\alpha , 0) \in A_\infty ] .$$
Observe that for $t < 1$,
\begin{eqnarray*}
M_t & = & m (A_t) + \one_{t < \tau} \int_{[0,1] \setminus A_t} 
   \P_{B_t} [(\alpha, 0) \in A_\infty ] \, d\alpha \\[2ex]
& = & \int_{l_0} f (B_{t \wedge \tau_\alpha} , \alpha ) d\alpha ,
\end{eqnarray*}
since the integrand in the last expression is 1 if $\alpha \in A_t$,
zero if $\alpha \notin A_t$ and $\tau \leq t$, and 
$\P_{B_t} [(\alpha , 0) \in A_\infty]$ otherwise.  

Lemma~\ref{lem integrate} is applicable to $M$ and $f$, 
with $(\S , \B , \mu) = (l_0 , \mbox{Borel}, m)$,
provided conditions $(i)$~-~$(iii)$ are satisfied.  
Checking these is not too hard.  Joint measurability, in fact
joint continuity, are clear, as is the vanishing of each
$f (\cdot , \alpha)$ on $\partial D$.  Condition~$(iii)$
is immediate from the Strong Markov Property.  By translation
invariance it suffices to check~$(ii)$ for $\alpha = 0$.  
This is Lemma~\ref{lem 3.45} below.
The conclusion of Lemma~\ref{lem integrate} is that $M_t$
is a continuous martingale with
\begin{equation} \label{eq 31}
<M>_1 = \int_0^1 \left | \int_0^1 \one_{s < \tau_\alpha}
   \grad f (B_s , \alpha)  \, d\alpha \right |^2 \; ds ,
\end{equation}
where $\tau_\alpha$ is the first time Brownian motion
hits either $\partial D$ or the ball of radius $\ee$
about the point $(\alpha , 0)$.  

The inequality~(\ref{eq to show}) will follow from 
equation~(\ref{eq 31}) along with some lemmas, stated
immediately below.

\begin{lem} \label{lem azuma}
For any continuous martingale $\{ M_t \}$ and any real $u$ and 
positive $t$ and $L$, 
$$\P [ M_t \geq u] \leq \exp (-(u-M_0)^2 / 2L) + \P [<M>_t > L] .$$
\end{lem}

\begin{lem} \label{pr 3}
$$d<M>_t \leq 4 \left | \int_0^\infty \grad f ((x , (B_t)_2) , 0)
   \, dx \right |^2 .$$
\end{lem}

\begin{lem} \label{lem 4}
There are constants $c_5, c_6$ and $c_7$ independent of $\ee$ for which
$$ (i) \;\;\; f ((0 , y) , 0) \leq c_5 {|\log |y|| \over |\log \ee|} ,$$
$$(ii) \;\;\; f ((|\alpha| , 0) , 0) \leq c_6 {e^{-\alpha} (1 + 
   (\log (1 / \alpha))^+ ) \over |\log \ee|} , $$
and
$$(iii) \;\;\; \int_{-\infty}^\infty f((0 , y) , \alpha) \, d\alpha 
   \leq c_7 {1-y \over |\log \ee|} . $$
\end{lem}

\begin{lem} \label{lem 3.45}
Fix $\ee > 0$.
For $x \in \R$ and $y \in (-1,1)$, let $\phi (x,y) = \P_{(x,y)} 
[\sigma < \tau]$,
where $\sigma$ is the hitting time on the $\ee$-ball around the origin
and $\tau$ is the first time $|Y_t| = 1$.  Then $\phi$ is Lipschitz.
\end{lem}

\begin{lem} \label{lem QV}
There is a constant $c_0$ for which
$$d<M>_t \leq c_0 \left ( {1 + |\log |Y_t|| \over |\log \ee|}
   \right )^2 .$$
\end{lem}

\begin{lem} \label{lem levy}
Let $L (x,t)$ be local time for a Brownian motion 
and let $U = \sup_x L(x,1)$.  Then for some constant $c_8$,
$\P [ U \geq u] \leq c_8 \exp (-u^2 / 2 c_8)$.
\end{lem}

\noindent{\sc Proof of }(\ref{eq to show}):  First, 
by Lemma~\ref{lem 4}, 
$$ M_0 = \int_0^1 f((0 , 0) , \alpha) \, d\alpha \leq {c_7
   \over |\log \ee|} .$$
Next, following the notation of the Lemma~\ref{lem levy}, 
let $H$ be the event $\{ U \leq |\log \ee| \}$.  According to
Lemma~\ref{lem levy} there is a constant $c_8$ for which
$$\P [H^c] \leq c_8 \exp (- (\log \ee)^2 / 2 c_8) .$$
Use Lemma~\ref{lem QV} to get
$$ <M>_1 \leq c_0 \left [ |\log \ee|^{-2} + \int_0^1 
   {\log^2 |Y_t| \over |\log \ee|^2 } \, dt \right ] = 
   c_0 |\log \ee|^{-2} \left [ 1 + \int_{-1}^1 \LL (y,1)
   \log^2 |y| \, dy \right ] $$
where $\LL$ is the local time function for the $y$-coordinate
of the Brownian motion.  On the event $H$, this is maximized
when $\LL(y,1) = |\log \ee|$ on the interval $[-|\log \ee| / 2 ,
|\log \ee| / 2]$ and has value at most 
$$ c_0 |\log \ee|^{-2} [2 (\log |2 \log \ee| )^2] .$$
Setting 
$$L = c_0 |\log \ee|^{-2} [2 (\log |2 \log \ee| )^2] $$
thus gives $\P [<M>_1 > L] \leq c_8 \exp (- (\log \ee)^2 / 2 c_8)$.
Plugging this value of $L$ into Lemma~\ref{lem azuma} yields
$$\P [M_1 \geq \theta] \leq \exp \left ( -{ (\theta - K |\log \ee|^{-1})^2
   \over 2L} \right ) + \P [ <M>_1 > L$$
where $K |\log \ee|^{-1}$ is an elementary estimate for the expected
measure of $l_0$ covered by the $\ee$-sausage to time 1.  This
is good enough to prove~(\ref{eq to show}) when $|\log \ee|$ is 
small compared to $\theta$.    $\Cox$

It remains to prove the six lemmas.

\noindent{\sc Proof of} Lemma~\ref{lem azuma}:  Routine; see 
the appendix.  $\Cox$

\noindent{\sc Proof of} Lemma~\ref{pr 3}:  First, it 
is evident that $f ((x,y),\alpha)$ depends only on 
$|x - \alpha|$ and $|y|$.  I claim that $f$ is decreasing in
$|x - \alpha|$ and $|y|$.  This may be verified
by coupling.  Let $\alpha \leq x \leq x'$ and $0 \leq
y \leq y'$ and couple two Brownian motions beginning at 
$(x , y)$ and $(x' , y')$ by letting the $x$-increments
be opposite until the two $x$-coordinates meet, then identical,
and the same for the $y$-coordinates.  Formally, if $W_t = (X_t , Y_t)$
is a Brownian motion, let $B_t = (x , y) + \int_0^t dW_s$
and $B_t' = (x' , y') + \int_0^t [ (1 - 2 \cdot \one_{s < \tau_x})
dX_s +  (1 - 2 \cdot \one_{s < \tau_y}) dY_s$, where
$\tau_x = \inf \{ t : (B_t)_1 \geq (x' + x) / 2 \}$ and
$\tau_y = \inf \{ t : (B_t)_2 \geq (y' + y) / 2 \}$.  This 
coupling has the property that $|B_t - (\alpha , 0)| \leq 
|B_t' - (\alpha , 0)|$ for all $t$, thus showing that
$f((x' , y') , \alpha) \leq f((x,y) ,\alpha)$.  

Use this claim to see that the signs of ${\partial f \over 
\partial x}$ and ${\partial f \over \partial y}$
are constant as $\alpha$ varies over $(-\infty , x)$, 
and are also constant when $\alpha$ varies over $(x , \infty)$.  
Hence
\begin{eqnarray*}
&& \left | \int \one_{t < \tau_\alpha} \grad f ((x,y) , \alpha) \, 
   d\mu (\alpha) \right | \\[2ex]
& \leq & \left | \int_{\alpha \leq x} \one_{t < \tau_\alpha} 
   \grad f ((x,y) , \alpha) \, d\mu (\alpha) \right | +
   \left | \int_{\alpha \geq x} \one_{t < \tau_\alpha} 
   \grad f ((x,y) , \alpha) \, d\mu (\alpha) \right |  \\[2ex]
& \leq & \left | \int_{\alpha \leq x} 
   \grad f ((x,y) , \alpha) \, d\mu (\alpha) \right | +
   \left | \int_{\alpha \geq x} 
   \grad f ((x,y) , \alpha) \, d\mu (\alpha) \right |  \\[2ex]
& \leq & 2 \left | \int_0^\infty \grad f ((x , y), 0) \, dx
   \right |
\end{eqnarray*}
by translation invariance.  Plugging this into the conclusion
of Lemma~\ref{lem integrate} gives the desired result.  $\Cox$

\noindent{\sc Proofs of Lemma}~\ref{lem 4} {\sc and}~\ref{lem 3.45}: 
Routine; see the appendix.  $\Cox$

\noindent{\sc Proof of Lemma }\ref{lem QV}:  Lemma~\ref{pr 3},
which gives estimates in terms of the gradient 
of $f$, may be turned into an estimate in terms of $f$ 
itself as follows.  First break down into the two coordinates:
$$d<M>_t \leq 4 |\int_0^\infty {\partial f \over \partial x}
   ((x , y) , 0) \, dx |^2 +
   4 |\int_0^\infty {\partial f \over \partial y}
   ((x , y) , 0) \, dx |^2 ,$$
where $y = Y_t$.  From the fundamental theorem of calculus, 
the square root of the first component is
\begin{equation} \label{eq 8}
2 |\int_0^\infty {\partial f \over \partial x}
   ((x , y) , 0) \, dx | = 2 (f((0 , y) , 0) - f((\infty ,
   y) , 0)) = 2 f ((0 , y) , 0) .
\end{equation}
To estimate the second component, assume without loss of generality
that $y > 0$ and write
\begin{equation} \label{eq 8.5}
2 |\int_0^\infty {\partial f \over \partial y}
   ((x , y) , 0) \, dx | = - {\partial \over \partial y} 
   \int_{-\infty}^\infty f((0 , y) , \alpha) \, d\alpha .
\end{equation}
For $\ee \in (0 , 1 - y)$, run a Brownian motion from 
$(0 , y + \ee)$ until it hits one of the horizontal lines 
$\R \times \{ 1 \}$ or $\R \times \{ y \}$.  The probability 
it hits $\R \times \{ y \}$ first is $1 - \ee / (1 - y)$.  
Sending $\ee$ to zero and using the Strong Markov Property gives
\begin{equation} \label{eq 9}
 - {\partial \over \partial y} 
   \int_{-\infty}^\infty f((0 , y) , \alpha) \, d\alpha =
(1 - y)^{-1} \int_{-\infty}^\infty f ((0 , y) , \alpha) 
   \, d\alpha .
\end{equation}
Combining~(\ref{eq 9}) and~(\ref{eq 8}) and
using Lemma~\ref{lem 4} gives $d<M>_t \leq 4 c_5^2 
(\log^2 |y| / \log^2 \ee) + c_7^2 / \log^2 \ee$, proving the
lemma with $c_0 = 4 c_5^2 + c_7^2$.  
$\Cox$

\noindent{\sc Proof of Lemma}~\ref{lem levy}: 
This is well known.  For example, it is stated and proved as Lemma~1
of~\cite{Cs}.    $\Cox$

\section{Proof of lower bounds}

In this section,~(\ref{eq lower}) and~(\ref{eq lomeasure})
are proved simultaneously.  For~(\ref{eq lomeasure}) I assume
that $\theta \leq 1 - 3\ee$, since otherwise this case is 
covered by~(\ref{eq lower}).  Let $N$ be a integer parameter to be 
named later.  Define a sequence of balls $\{ \C_j : 1 \leq j \leq 
N^2 \}$, each of radius $1/N$ and centered on the $x$-axis with 
ordinates
$$1/N , 2/N , \ldots , (N-1)/N , 1 , (N-1)/N , \ldots , 1/N , 0 , 1/N , 
   \ldots , 1 , \ldots , 1/N , 0 , \ldots $$
running back and forth along the unit interval a total of $N$ times.
Let $H_N$ be the event that 
$$B_{j/N^2} \in \C_j \mbox{ for all } j \leq N^2 \; . $$
Let $c_9 = \log \P_{(-1/N , 0)} [ B_{1/N^2} \in \C_1 ]$.  By
the Markov property, 
$$ \P [B_{j / N^2} \in \C_j \| B_{i/N^2} \in \C_i : 1 \leq i < j ]
   \geq \min_{z \in \C_{j-1}} \P_z [ B_{1/N^2} \in \C_j ] ,$$
which is easily seen to be minimized when $j=1$ and $z = (-1/N , 0)$, 
leading to 
\begin{equation} \label{eq 21}
\P [ H_N ] \geq \exp (-c_9 N^2 ) .
\end{equation}

Next, bound from below the conditional probabilities of the events 
$\{ R_\ee \supseteq [0,1] \}$ and of $\{ \Xi \geq \theta \}$ given $H_N$
(recall $\Xi$ is the measure of $R_\ee \cap ([0,1] \times \{ 0 \} )$).
Let $b = \lceil 1/\ee \rceil$ be the least integer exceeding $1/\ee$.
For $1 \leq j \leq b$, let $z_j = (j \ee , 0)$.
If $z_j \in R_{\ee/2}$ then $[z_j - \ee/2 , 1 \wedge (z_j + \ee/2)]
\times \{ 0 \} \subseteq R_\ee$; thus $\Xi \geq \ee ( \# \{ j :
z_j \in R_{\ee / 2} \} - 1 )$ and if every $z_j$ is in $R_{\ee / 2}$
then $\{ R_\ee \supseteq [0,1] \}$.  The following routine lemma is
proved in the appendix.

\begin{lem} \label{lem bridge}
There is a universal constant $K$ such that
if $q_1, q_2$ and $q_3$ are points in $\R^2$ at pairwise distances of no 
more than $3  / N$ and $\{ W_t : 0 \leq t \leq 1/N^2 \}$ is a Brownian
bridge with $W_0 = q_1$ and $W_{1/N^2} = q_2$, then
$$\P [ |W_s - q_3| < \dd \mbox{ for some } s \leq 1 / N^2 ] \geq
	{1 \over K |\log \dd| } $$
for all $N \geq 1$ and $\dd < 1/2$.
$\Cox$
\end{lem}

Continuing the proof of the lower bounds, fix $j$ and condition
on the $\sf$ $\sigma ( B_{i/N^2} : 1 \leq i \leq N^2 )$.
On the event $H_N$ there are $N$ values of $i$ for which the 
pairwise distances between $B_{i/N^2} , B_{(i+1)/N^2}$ and $z_j$
are all at most $3/N$.  It follows from Lemma~\ref{lem bridge}
that on $H_N$, 
$$\P [ z_j \notin R_{\ee/2} ] \leq \left ( 1 - {1 \over K 
   |\log (\ee/2)|} \right )^N .$$
Choosing $N = \lceil \gamma K |\log (\ee / 2)| \rceil$ now gives 
$$\P [ z_j \notin R_{\ee/2} \| H_N ] \leq e^{-\gamma} . $$
When $e^{-\gamma} \leq (1 - \theta) / 3$ this gives
$\E \left [ \# \{ j < b : z_j \notin R_{\ee / 2} \} \| H_N \right ]
\leq (1 - \theta) (b-1) / 3$ and hence
$$\P [ \# \{ j < b : z_j \notin R_{\ee / 2} \} \geq (1 - \theta) b 
   - 1 \| H_N ] \leq {(1 - \theta) (b-1) / 3 \over (1 - \theta) b - 1}$$
which is at most $1/2$ when $\theta \leq 1 - 3\ee$.  Thus
choosing $\gamma = |\log (1-\theta)/3|$ gives
$$\P [ \Xi \geq \theta ] \geq \P  [ \# \{ j < b : z_j \notin 
   R_{\ee / 2} \} < (1 - \theta) b - 1 ] \geq {1 \over 2} 
   \P [ H_N ] \geq {1 \over 2} \exp (- c_9 \gamma^2 K^2 
   |\log \ee / 2|^2) .$$
This is at least $\exp ( - c_4 |\log \ee|^2 |\log (1 - \theta)/3|^2)$
when $c_4 = 2 c_9 K^2$ and $\ee$ is small enough, proving~(\ref{eq lomeasure}).

Choosing instead $N = \lceil \gamma K |\log (\ee / 2)| |\log \ee| \rceil$ 
makes
$$\P [ z_j \notin R_{\ee/2} \| H_N ] \leq \ee^\gamma $$
and for $\gamma = 2$ this makes $\P [R_\ee \supseteq [0,1] 
\times \{ 0 \} \| H_N ]$ at least $1 - 
\lceil \ee^{-1} \rceil \ee^2 \geq exp (- c |log \ee|^4)$ for any $c$ and
small $\ee$.
For such a $\gamma$,~(\ref{eq 21}) gives 
\begin{eqnarray*}
\P [H_N] & \geq & \exp (- c_9 (2 K |\log \ee| )^4) \\[2ex]
& \geq & \exp (- (1/2) c_4 |\log \ee|^4)
\end{eqnarray*}
for $c_4 = 32 K^4 c_9$ and $\ee$ small enough; thus $\P [ R_\ee \supseteq [0,1]
\times \{ 0 \}] \geq \exp (- c_4 |\log \ee|^4)$, proving~(\ref{eq lower}).
$\Cox$

\section{Appendix: proofs of easy propositions} 

\noindent{\sc Proof of Proposition}~\ref{pr COV}:  Continuity is 
clear, and the vanishing
of the Laplacian of $f$ shows that each $M^{(\alpha)}$ is 
a local martingale, hence a martingale since it is bounded.
In fact the formula for the brackets is clear too, but
to be pedantic we substitute for $f$ a sequence of functions
genuinely in $C^2$ and apply It\^o's  formula.  
For $\ee > 0$ let $D_{\alpha , \ee} = 
\{x \in D : f(x , \alpha) \in [\ee , 1-\ee] \}$.  Use a
partition of unity to get functions $f^{\alpha, \ee} : \R^2 
\rightarrow [0,1]$ with continuous, uniformly bounded second
derivatives and such that $f^{\alpha , \ee} = f$
on $D_{\alpha , \ee}$.  Let $\tau_{\alpha , \ee} = \inf
\{ t : B_t \notin D_{\alpha , \ee} \}$.  The martingale
$$M_t^{\alpha , \ee} := f (B_{t \wedge \tau_{\alpha , \ee}} ,
   \alpha)$$
is a continuous martingale and it is clear from It\^o's formula
that the covariances satisfy
\begin{equation} \label{eq 2}
\E \left [ (M_t^{\alpha , \ee} - M_s^{\alpha , \ee}) (M_t^{\beta 
   , \ee} - M_s^{\beta , \ee}) \| \F_s \right ]
   = \int_s^t 
   \one_{ r < \tau_{\alpha , \ee}} \one_{ r < \tau_{\beta , \ee}}
   \grad f(B_r , \alpha) \cdot \grad f(B_r , \beta) \, dr .
\end{equation}
Now fix $\alpha$ and $\beta$ and let $\ee$ go to zero.  Since
$\tau_{\alpha , \ee} \uparrow \tau_\alpha$ almost surely
it follows that $M_r^{\alpha , \ee} \rightarrow M_r^{(\alpha)}$ 
almost surely and $\one_{r < \tau_{\alpha , \ee}} \uparrow 
\one_{r < \tau_\alpha}$ almost surely.  The same is true 
with $\beta$ in place of $\alpha$.  The expression inside the
expectation (resp. integral) on the left (resp. right) side 
of~(\ref{eq 2}) is bounded, hence sending $\ee$ to zero 
simply removes it typographically from both sides, and
establishes the lemma.   $\Cox$

\noindent{\sc Proof of Proposition}~\ref{pr together}:  The 
integral of bounded continuous functions
is continous.  To check that $M$ is a martingale, first observe
that 
$$\int \E (|M_t^{(\alpha)} - M_s^{(\alpha)}| \| \F_s ) \, d\mu (\alpha) 
   \leq \mu (\S) \sqrt{C' (t-s)} $$
and hence one may interchange an integral against $\mu$ with
an expectation given $\F_s$ to get
$$\E [\int (M_t^{(\alpha)} - M_s^{(\alpha)}) \, d\mu (\alpha) \| \F_s] 
   = \int \E (M_t^{(\alpha)} - M_s^{(\alpha)} \| \F_s) \, d\mu (\alpha) .$$ 
Thus,
$$ \E (M_t - M_s \| \F_s) = 0.$$
To compute the quadratic variation, first observe that
$$\E \int \int \E \left [ |(M_t^{(\alpha)} - M_s^{(\alpha)})
   (M_t^{(\beta)} - M_s^{(\beta)})| \| \F_s \right ] \, 
   d\mu (\alpha) \, d\mu (\beta) \leq 1 .$$
This justifies the third of the following equalities, the second
being justified by absolute integrability:
\begin{eqnarray*}
[ (M_t - M_s)^2 \| \F_s ] & = & \E \left [ \left ( \int
   (M_t^{(\alpha)} - M_s^{(\alpha)}) \, d\mu (\alpha) \right )^2
   \| \F_s \right ] \\[2ex]
& = & \E \left [ {\int \int}_{\S^2} (M_t^{(\alpha)} - M_s^{(\alpha)})
   (M_t^{(\beta)} - M_s^{(\beta)}) \, d\mu (\alpha) \, d\mu (\beta)
   \| \F_s \right ] \\[2ex]
& = & {\int \int}_{\S^2} \E \left ( (M_t^{(\alpha)} - M_s^{(\alpha)})
   (M_t^{(\beta)} - M_s^{(\beta)}) \| \F_s \right ) 
   \, d\mu (\alpha) \, d\mu (\beta) \\[2ex]
& = & {\int \int}_{\S^2} \left ( [M^{(\alpha)} , M^{(\beta)}]_t - 
   [M^{(\alpha)} , M^{(\beta)}]_s \right ) \, d\mu (\alpha) \, 
   d\mu (\beta) \\[2ex]
& = & \int_s^t \left [ {\int \int}_{\S^2} d[M^{(\alpha)} , M^{(\beta)} 
   ]_r \, d\mu(\alpha) \, d\mu (\beta) \right ] 
\end{eqnarray*}
again by absolute integrability.  $\Cox$

\noindent{\sc Proof of Lemma }\ref{lem azuma}: Let $\tau = \inf \{ s :
<M>_s \geq L \}$.  Then
$$\exp ( \lambda M_{s \wedge \tau} - {1\over 2} \lambda^2 
   <M>_{s \wedge \tau} )$$
is a martingale and hence 
\begin{eqnarray*}
\E \exp (\lambda (M_{t \wedge \tau} - M_0)) & \leq & \E \exp \left (
   \lambda (M_{t \wedge \tau} - M_0) - {\lambda^2 \over 2} 
   (<M>_{t \wedge \tau} - L) \right ) \\[3ex]
& \leq & \exp ({\lambda^2 \over 2} L) .
\end{eqnarray*}
By Markov's inequality, 
\begin{eqnarray*}
\P [M_{t \wedge \tau} - M_0 \geq u - M_0] & \leq & \exp (- \lambda (u - M_0))
  \E \exp (\lambda (M_{t \wedge \tau} - M_0)) \\[2ex]
& \leq & \exp ({\lambda^2 \over 2} L - \lambda (u - M_0)) .
\end{eqnarray*}
Plugging in $\lambda = (u - M_0) / L$ gives
$$\P [M_{t \wedge \tau} \geq u] \leq \exp (-(u - M_0)^2 / 2L) $$
which proves the lemma since $\{ M_t \geq u \} \subseteq \{ M_{t \wedge
\tau} \geq u \} \cup \{ <M>_t > L \}$.   $\Cox$

\noindent{\sc Proof of Lemma }\ref{lem 4}:
Assume without loss of generality that 
$\ee < 1/4$.  Let $\C_\ee$ be the circle of radius $\ee$ centered
at the origin, and $\tau_\ee$ the hitting time of this circle.
These are useful stopping times, since for $a < |z| < b$,
$$ \P_z [\tau_a < \tau_b] = {\log (b/|z|) \over \log (b/a)} .$$
The minimal value $\P_z [\tau < \tau_{\C_{1/2}}]$ as $z$ ranges
over points of modulus $3/4$ is some constant $p_0$.  Given
a Brownian motion, construct a continuous time
non-Markovian birth and death process $\{ X_t \}$ taking values 
in $\{ 0,1,2,3 \}$ by letting $X_t$ be $0,1,2$ or $3$ according to 
whether Brownian motion has most recently hit $\C_\ee , \C_{1/2} , 
\C_{3/4}$ or $\partial D$ first respectively.  Let $\{ X_n \}$
be the discrete time birth and death process consisting of 
successive values of $\{ X_t \}$.  Let $\F_n = \sigma (X_1 , 
\ldots , X_n)$.  Then 
\begin{eqnarray*}
\P (X_{n+1} = 0 \| \F_n , X_n = 1) & = & {\log (3/2) \over 
   \log (3/4\ee)} \, \mbox{ and} \\[2ex] 
\P (X_{n+1} = 3 \| \F_n , X_n = 2) & \geq & p_0 .
\end{eqnarray*}
Couple this to a Markov birth and death chain whose transition
probability from 2 to 3 is precisely $p_0$, to see that
for this chain, the probability of hitting 0 before 3 given
$X_1 = 1$ is at most $K / |\log \ee|$ for some constant $K$.
This implies that for $y > 1/2$, 
$$f ((0,y) ,0) \leq \P_{(0,y)} [\tau_{1/2} < \tau ] { K
   \over |\log \ee|} \, .$$
Since $\P_{(0,y)} [\tau_{1/2} < \tau ] < 2(1-y) < |\log y|$,
this verifies~$(i)$ in the case $y > 1/2$.  For 
$y \leq 1/2$, use 
$$f ((0,y) ,0) \leq \P_{(0,y)} [\tau_{1/2} < \tau_\ee ] { K
   \over |\log \ee|} + \P_{(0,y)} [\tau_\ee < \tau_{1/2}] $$
to verify~$(i)$ in the remaining cases.  

This also verifies~$(ii)$
in the case $\alpha < 1/2$.  For $1/2 \leq \alpha \leq 2$, 
the numerator on the RHS of~$(ii)$ is bounded away from 0 and the 
LHS is decreasing in $\alpha$, so $c_6$ may be chosen a little 
larger to make~$(ii)$ hold.  Finally, for $\alpha > 2$, note that 
$$\P_{(k+1 , y)} [\sigma_k < \tau ] < 1/3 < e^{-1} ,$$
for any $y \in [0,1]$, where $\sigma_k$ is the hitting time on 
the vertical line $x = k$.  Applying this successively gives 
the factor $e^{-\alpha}$ and verifies~$(ii)$.  

To see~$(iii)$ for $y \leq \ee$, use monotonicity to get
$$f ((0 , \ee) , \alpha) \leq f((0 , 0) , \alpha) = f((\alpha , 0) , 0)$$
and integrate~$(ii)$.  Finally, for $\ee < y < 1$, let
$T$ be the hitting time on the horizontal line $y=\ee$ and observe that
$$ \int_{-\infty}^\infty f((0 , y) , \alpha) \, d\alpha 
   = \P [T < \tau ] 
   \int_{-\infty}^\infty f((0 , \ee) , \alpha) \, d\alpha .$$
$\Cox$

\noindent{\sc Proof of Lemma}~\ref{lem 3.45}:  First let $(x,y)$ be a point
inside the strip $|y| \leq 1-\ee$ but outside the ball $x^2 + y^2 < 4 \ee^2$.
Let $V$ be the circle of radius $\ee/2$ centered at $(x,y)$.  If
$(x' , y')$ is another point inside $V$ and $\mu_{x' ,y' }$ is 
the hitting distribution on $V$ starting from $(x' , y')$, then
the total variation distance from $\mu_{x' , y'}$ to uniform is 
$O( \ee^{-1} |(x,y) - (x' , y')|)$ as $(x' , y') \rightarrow (x,y)$.  
The strong Markov property shows that $|\phi (x' , y') - \phi (x,y)|$
is bounded by this total variation distance; the truth of this
statement in the limit $(x' , y') \rightarrow (x,y)$ implies 
that $\phi $ is Lipschitz with a possibly greater constant in a neighborhood
of $(x,y)$.  Inverting in the circle $x_2 + y^2 = \ee^2$ and using a 
similar argument shows that $\phi$ is also Lipschitz in a neighborhood
of $\infty$ on the complex sphere $C^*$, so compactness of the region 
$\{ |y| \leq 1 - \ee \} \setminus \{ x^2 + y^2 < 4 \ee^2 \}$ in $C^*$
shows that $\phi$ is Lipschitz in that region. 

The two other cases are also easy.  If $\ee \leq |(x,y)| \leq 2 \ee$ then
let $\tau$ be the first time that $|(X_y , Y_t)| = \ee$ or $2 \ee$.
Let $\mu_{x,y}$ be the subprobability hitting measure on the
outer circle $V \deq \{ (u,v) : u^2 + v^2 = 4\ee^2 \}$, i.e., $\mu (B)
= \P[ (X_\tau , Y_\tau ) \in B \cap V]$.  
If $(x' , y')$ is another point of modulus between $\ee$ and $2 \ee$, 
then $|\phi (x' , y') - \phi (x,y)|$ is at most the total variation distance
between $\mu_{x,y}$ and $\mu_{x' , y'}$.  Again, since $\ee$
is fixed, an easy computation shows this to be $O(|(x,y) - (x' , y')|)$.
The case where $y > 1 - \ee$ is handled the same way, except there
the outer circle is replaced by the line $y = 1 - 2 \ee$.   $\Cox$

\noindent{\sc Proof of Lemma}~\ref{lem bridge}:  Assume without loss
of generality that $q_3$ is the origin and dilate time by $N^2$ and 
space by $N$, so that we have a Brownian bridge $\{ X_t : 0 \leq t 
\leq 1 \}$ started at $p_1$ and stopped at $p_2$, with $|p_1| , |p_2|
\leq 3$.  We use the following change of measure formula for the law 
of this Brownian bridge.  Let $\tau \leq 1$ be any stopping time, let
$\P_x$ denote the law of Brownian motion from $x$ and $\Q$ denote
the law of the bridge from $q_1$ at time 0 to $q_2$ at time 1.
Then for any event $G \in \F_\tau$, 
$$ \Q (G) = {\int \one_G \exp (|q_2 - B_\tau|^2 / 2 (1-\tau)) 
   \, d\P_{q_1} \over \exp (|q_2 - q_1|^2 / 2) } .$$

The probability of the bridge entering the dilated 
$N \dd$-ball around the origin is of course at least the probability
it enters the $\dd$-ball around the origin before time $3/4$.
If $\{ B_t \}$ is an unconditioned Brownian motion 
starting from $p_1$ and $\tau$ is the first time $\{ B_t \}$
enters the $\dd$-ball around the origin, then $\P [\tau \leq 3/4]$
is at least $K_1 / |\log \dd|$ for a universal constant $C_1$.
Conditioned on $\tau$ and $B_\tau$, the density of $B_1$ at $q_2$
is at least $(2\pi)^{-1} e^{-18}$.  Plugging into the change of measure
formula proves the lemma.  $\Cox$


\begin{thebibliography}{YMN}

\bibitem{BPP}
Benjamini, I., Pemantle, R. and Peres, Y. (1994).  Martin capacity for
Markov chains.  {\em Ann. Probab., to appear.}

\bibitem{Bu3} 
Burdzy, K. (1994).  Labyrinth dimension of Brownian trace.  {\em
Preprint}.

\bibitem{BL}
Burdzy and Lawler (1990). Nonintersection exponents for Brownian 
paths II: estimates and application to a random fractal.  
{\em Ann. Probab.} {\bf 18} 981 - 1009.

\bibitem{Cs}
Cs\'aki, E. (1989).  An integral test for the supremum of Wiener 
local time.  {\em Prob. Th. Rel. Fields} {\bf 83} 207 - 217.

\bibitem{DEK}
Dvoretzky, A., Erd\"os, P. and Kakutani, S. (1950).  Double points of paths 
of Brownian motion in $n$-space.  {\em Acta Sci. Math.} {\bf 12} 75 - 81.

\bibitem{FS}
Fitzsimmons, P.J.  and Salisbury, T. (1989). 
Capacity and energy for multiparameter Markov processes.
{\em Ann.\ Inst.\ Henri Poincar\`e, Probab.\/} {\bf 25} 325-350.

\bibitem{Kak}
Kakutani, S. (1944).  Two dimensional Brownian motion and harmonic functions.
  {\em Proc. Imp. Acad. Tokyo} {\bf 20}, 648-652.

\bibitem{La}
Lawler, G. (1992).  On the covering time of a disc by simple random
walk in two dimensions.  In: Seminar on stochastic processes, 1992,
R. Bass and K. Burdzy managing editors.  Birkh\"auser: Boston.

\bibitem{LeG}
Le Gall, J. F. (1991).  Some properties of planar Brownian motion.
{\em Springer Lecture Notes in Mathematics} {\bf 1527} 112 - 234.

\bibitem{Lev}
L\'evy, P. (1940).  Le mouvement brownien plan.  {\em Amer. J. Math.}
{\bf 62} 487 - 550.

\bibitem{Anon}
Meyre, T. and Werner, W. (1994).  Estimation asymptotique du rayon du 
plus grand disque couvert par la saucisse de Wiener plane.  
Stochastics and Stochastics Reports {\bf 48}, 45 - 59.

\bibitem{Ra}
Ray, D. (1963).  Sojourn times and the exact Hausdorff measure of the sample 
path for planar Brownian motion.  {\em Trans. AMS} {\bf 106} 436 - 444.

\bibitem{Re}
R\'ev\'esz, P. (1990).  Estimates on the largest disc covered by a 
random walk.  {\em Ann. Probab.} {\bf 18} 1784- 1789.

\bibitem{Sa}
Salisbury, T. (1994).  Energy, and intersections of Markov chains.
In: IMA volume on Random discrete structures, D. Aldous and R. Pemantle Eds.,
to appear.

\end{thebibliography}
\end{document}